\documentclass{commat}

\usepackage{graphicx}

\def\erf{\mathop{\mathrm{erf}}}
\newcommand{\fact}[1]{#1!}

\allowdisplaybreaks

\title{Matrix continued fractions and Expansions of the Error Function}

\author{S. Mennou, A. Chillali and A. Kacha}

\affiliation{
    \address{
        S. Mennou, A. Kacha --
        Ibn Toufail University, Science faculty, Laboratory AMGNCA, Kenitra 14000, Morocco.
        }
    \email{saidmennou@yahoo.fr}
    \email{ali.kacha@yahoo.fr}
    \address{
        A. Chillali --
        Sidi Mohamed Ben Abdellah University, Polydisciplinary Faculty , LSI, Taza, Morocco.
        }
    \email{abdelhakim.chillali@usmba.ac.ma}
    }

\abstract{
In this paper we recall some results and some criteria on the
convergence of matrix continued fractions. The aim of this paper is
to give some properties and results of continued fractions with
matrix arguments. Then we give continued fraction expansions of the error
function $\erf(A)$ where A is a matrix. At the end, some numerical
examples illustrating the theoretical results are discussed.	
}



\VOLUME{31}
\YEAR{2023}
\NUMBER{1}
\firstpage{257}
\DOI{https://doi.org/10.46298/cm.10395}

\begin{paper}

\section{Introduction and motivation}

The theory of continued fractions has been a~topic of great interest
over the last two hundred years. The basic idea of this theory over
real numbers is to give an approximation of various real numbers by the rationals.
A continued fraction is an expression obtained
through an iterative process of representing a~number as the sum of
its integer part and the reciprocal of another number, then writing
this other number as the sum of its integer part and another
reciprocal, and so on. One of the main reasons why continued
fractions are so useful in computations is that they often provide
representations for transcendental functions that are much more
generally valid than the classical representation by, say, the power
series. Further, in the convergent case, continued fraction expansions
 have the advantage that they converge more rapidly than
other numerical algorithms.

Recently, the extension of continued fraction theory from real
numbers to the matrix case has seen several developments and
interesting applications (see \cite{1}, \cite{3}, \cite{6}). Since calculations
involving matrix valued functions with matrix arguments are feasible
with large computers, it will be an interesting attempt to develop such a~matrix
 theory. In this direction, and generally in a~Banach space,
few convergence results on non-commutative continued fraction are
known.

Two theorems are stated in~\cite{10}, where Wynn reviews many aspects
of the theory of continued fractions, whose elements do not commute
under a~multiplication law. In Banach space, extensions of
Worpitsky's have been proven by Haydan~\cite{2} and Negoescu~\cite{8}.

In \cite{9}, the authors give several convergence criteria on non-commutative
continued fractions whose arguments are $m \times m$ matrices of the
form $K(B_{n}/A_{n})$. 

The error function $\erf$  is a~special function that is important
since it appears in the solutions of diffusion problems in heat,
mass and momentum transfer, probability theory, theory of errors and
various branches of mathematical physics. The closely related
Fresnel integrals, which are fundamental in the theory of optics,
can be derived directly from the error function.

\section{Preliminaries and notations}
Throughout this paper, we denote by $\mathcal{M}_{m}$  the set of $m\times m$ 
 real (or complex) matrices endowed with the subordinate matrix
infinity norm defined by,
\begin{align*}
\forall A = (a_{i,j}),\ A \in \mathcal{M}_{m},\ \lVert A\rVert
= \max_{1\leq i\leq m}\sum_{j = 1}^{m} |a_{i,j}|.
\end{align*}
This norm satisfies the inequality
$$
\lVert AB\rVert \leq \lVert A\rVert \lVert B\rVert\,.
$$

Let $A\in \mathcal{M}_{m}$, $A$  is said to be positive semi-definite
(resp. positive definite) if $A$  is symmetric and
$$
\forall x\in \mathbb{R}^{m}, (Ax, x)\geq 0 \ (\text{resp. } \forall
x\in \mathbb{R}^{m},\ x\neq 0,\ (Ax, x)> 0)
$$
where $(\cdot,\cdot)$  denotes the standard scalar product of
$\mathbb{R}^{m}$  defined by

$$\forall x = (x_{1}, \ldots, x_{m})\in \mathbb{R}^{m}, \forall
y = (y_{1},\ldots, y_{m}) \in \mathbb{R}^{m}: (x, y) = \sum_{i = 1}^{m}
x_{i}y_{i}\,.$$

For any $A,B \in \mathcal{M}_{m}$  with $B$  invertible, we write $A/B:= B^{-1}A$, in particular,
if $A = I$, where $I$  is the $m^{th}$  order identity matrix, then we
write $I/B = B^{-1}$. It is clear that for any invertible matrix $C$,
 we have
$$\frac{CA}{CB} = \frac{A}{B}\,.$$

\begin{definition}
Let $(A_{n})_{n\geq 0}, (B_{n})_{n\geq 1}$  be two nonzero sequences
of $\mathcal{M}_{m}$. The continued fraction of $(A_{n})$  and
$(B_{n})$, denoted by $K(B_{n}/A_{n})$, is the quantity

\begin{equation*}
A_{0}+\frac{B_{1}}{A_{1}+\frac{B_{2}}{A_{2+\ldots}}} = \Bigl[A_{0};\frac{B_{1}}{A_{1}},\frac{B_{2}}{A_{2}},\ldots\Bigr].
\end{equation*}

Sometimes, we use the notation
$\left[A_{0};\frac{B_{k}}{A_{k}}\right]_{k = 1}^{+\infty}$  or
$K(B_{n}/A_{n})$, where
\begin{equation*}
\left[A_{0};\frac{B_{k}}{A_{k}}\right]_{k = 1}^{n} = \left[A_{0};\frac{B_{1}}{A_{1}},
\frac{B_{2}}{A_{2}},\ldots, \frac{B_{n}}{A_{n}}\right].
\end{equation*}

The fractions $\frac{B_{n}}{A_{n}}$  and
$\frac{P_{n}}{Q_{n}}:= \left[A_{0};\frac{B_{k}}{A_{k}}\right]_{k = 1}^{n}$ 
are called, respectively, the $n^{th}$  partial quotient and the
$n^{th}$ 
convergent of the continued fraction $K(B_{n}/A_{n})$.

When $B_{n} = I$  for all $n\geq 1$, then $K(I/A_{n})$ is called an ordinary continued fraction.
The following proposition gives an
adequate method to calculate $K(B_{n}/A_{n})$.
\end{definition}

\begin{proposition}
The elements $(P_{n})_{n\geq -1}$  and $(Q_{n})_{n\geq -1}$ of the
$n^{th}$  convergent of $K(B_{n}/A_{n})$  are given by the
relationships
$$
\begin{cases}
 P_{-1} =  I,  &P_{0} = A_{0}\\
 Q_{-1} =  0, &Q_{0} = I
 \end{cases}
 \text{ and }
 \begin{cases}
 P_{n}  =  A_{n}P_{n-1}+B_{n}P_{n-2} \\
 Q_{n}  =  A_{n}Q_{n-1}+B_{n}Q_{n-2}
 \end{cases}\,,
n\geq 1\,.
$$
\end{proposition}

\begin{proof}
This can be done by induction.
\end{proof}

 The proof of the next Proposition is elementary and we leave it to the reader.
 \begin{proposition}
 For any two matrices $C$  and $D$  with $C$  invertible, we
have
\begin{equation}
C\left[A_{0};\frac{B_{k}}{A_{k}}\right]_{k = 1}^{n}D = \left[CA_{0}D;\frac{B_{1}D}{A_{1}C^{-1}},\frac{B_{2}C^{-1}}{A_{2}},
\frac{B_{k}}{A_{k}}\right]_{k = 3}^{n}.
\end{equation}
\end{proposition}

The continued fraction $K(B_{n}/A_{n})$  converges in $\mathcal{M}_{m}$  if the sequence $$(F_{n}) =
\Bigl(\frac{P_{n}}{Qn}\Bigr) = (Q^{-1}_{n}P_{n})$$ converges in
$\mathcal{M}_{m}$  in the sense that there exists a~matrix $F \in
\mathcal{M}_{m}$  such that
\[
    \lim_{n \to +\infty}||F_{n}-F|| = 0.
\]
In the opposite case, we say that $K(B_{n}/A_{n})$  is divergent. It
is clear that
\begin{equation}
\frac{P_{n}}{Q_{n}} = A_{0}+
\sum\limits_{i=1}^n\Bigl(\frac{P_{i}}{Q_{i}}-\frac{P_{i-1}}{Q_{i-1}}\Bigr)\,,
\label{moneq}
\end{equation}
and from~\eqref{moneq}, we see that the continued fraction
$K(B_{n}/A_{n})$  converges in $\mathcal{M}_{m}$  if and only if the
series
$\sum\limits_{\substack{n = 1}}^{+\infty}{(\frac{P_{n}}{Q_{n}}-\frac{P_{n-1}}{Q_{n-1}})}$ 
converges in $\mathcal{M}_{m}$.
 
\begin{definition}
Let $(A_{n})$, $(B_{n})$, $(C_{n})$  and $(D_{n})$  be four sequences
of matrices. We say that the continued fractions $K(B_{n}/A_{n})$ 
and $K(D_{n}/C_{n})$  are equivalent if we have $F_{n} = G_{n}$  for
all $n\geq1$, where $F_{n}$  and $G_{n}$  are the $n^{th}$  convergent
of $K(B_{n}/A_{n})$  and $K(D_{n}/C_{n})$  respectively.
\end{definition}

The following lemma characterizes equivalence of continued fractions.
\begin{lemma}[\cite{4}]
\label{lem} Let $(r_{n})$  be a~non-zero sequence of real numbers.
The continued fractions
\begin{equation*}
\left[a_{0};\frac{r_{1}b_{1}}{r_{1}a_{1}},
\frac{r_{2}r_{1}b_{2}}{r_{2}a_{2}}, \ldots,
\frac{r_{n}r_{n-1}b_{n}}{r_{n}a_{n}}, \ldots \right]\text{ and }
\left[a_{0};\frac{b_{1}}{a_{1}}, \frac{b_{2}}{a_{2}}, \ldots,
\frac{b_{n}}{a_{n}},\ldots \right]
 \end{equation*}
 are equivalent.
\end{lemma}

We also recall the following Lemma. From the expansion of a~function given by its Taylor
series, we give the expansion in continued fractions of the series that was established by Euler.
\begin{lemma}[\cite{5}] \label{lemma}
Let $f$  be a~function with Taylor series expansion $f(x) =
\sum\limits_{\substack{n = 0}}^{+\infty}{c_{n}x^{n}}$  in $D\subset R$.
Then, the expansion in continued fraction of $f(x)$  is
\begin{align*}
f(x) &= \left[ c_{0};\frac{c_{1}x}{1}, \frac{-c_{2}x}{c_{1} + c_{2}x}, \frac{-c_{1}c_{3}x}{c_{2} + c_{3}x},\ldots
 \frac{-c_{n-2}c_{n}x}{c_{n-1}+ c_{n}x},\ldots \right]\\
 &= \left[ \frac{c_{0}}{1}, \frac{-c_{1}x}{c_{0} + c_{1}x}, \frac{-c_{0}c_{2}x}{c_{1} + c_{2}x},
 \frac{-c_{1}c_{3}x}{c_{2} + c_{3}x},\ldots, \frac{-c_{n-2}c_{n}x}{c_{n-1}+ c_{n}x},\ldots \right].
\end{align*}
\end{lemma}

\begin{remark}
Let $(A_{n})$  and $(B_{n})$  be two sequences of $\mathcal{M}_{m}$. Then
we notice that we can write the first convergents of the continued
fraction
$K(B_{n}/A_{n})$  by:
$$
F_{1} = A_{0}+A^{-1}_{1}B_{1} = A_{0}+(B_{1}^{-1} A_{1})^{-1}.
$$
 $$
F_{2} = A_{0}+ (A_{1} + A^{-1}_{2}B_{2})^{-1}B_{1} \\
 = A_{0} +(B_{1}^{-1}A_{1}+(B_{2}^{-1}A_{2}B_{1})^{-1})^{-1}.
$$
 If we put, $A^{\ast}_{1} = B_{1}^{-1}A_{1}$  and $A^{\ast}_{2} = B_{2}^{-1}A_{2}B_{1}$, we have
 $$
F_{1} = A_{0} + \frac{I}{A_{1}^{\ast}}, F_{2} = A_{0} + \frac{I}{A^{\ast}_{1}+\frac{I}{A^{\ast}_{2}}}.
$$
 Generally, we prove by a~recurrence that if we put
for all $k \geq 1$, 
$$
A^{\ast}_{2k} = (B_{2k} \ldots B_{2})^{-1} A_{2k} B_{2k-1} \ldots B_{1}
$$
and
$$
A^{\ast}_{2k+1} = (B_{2k+1} \ldots B_{1})^{-1}
A_{2k+1} B_{2k} \ldots B_{2},
$$
then the continued fractions
$A_{0}+K(B_{n}/A_{n})$  and $A_{0}+K(I/A^{\ast}_{n})$  are equivalent.
\end{remark}

So, the convergence of one of these continued fractions implies the
convergence of the other continued fraction.

\begin{theorem}[\cite{12}]
\label{th} Let all the elements of $A_{n} (n = 1, 2, \ldots)$  be
positive, namely, $A_{n}$  are positive matrices for all $n$, then
the matrix continued fraction $K(I/A_{n})$  converges if and only if
$\sum\limits_{\substack{n = 0}}^{+\infty}{||A_{n}||} = \infty$.
\end{theorem}

We will use the following Theorem to prove our main result.

\begin{theorem}[\cite{9}]
\label{cov} Let $(A_{n})$, $(B_{n})$  be two sequences of
$\mathcal{M}_{m}$. If
$$
||(B_{2k-2} \ldots B_{2})^{-1} A^{-1}
_{2k-1}B_{2k-1} \ldots B_{1}|| \leq \alpha
$$
and
$$
|| (B_{2k-1} \ldots B_{1})^{-1} A^{-1} _{2k}
B_{2k} \ldots B_{2})||  \leq \beta
$$
 for all $k\geq 1$, where $0 <\alpha< 1, 0<\beta< 1 \; and \; \alpha\beta\leq1/4$, then the continued fraction $K(Bn/An)$ 
converges in $\mathcal{M}_{m}$.
\end{theorem}

 We need to present the following Proposition:

 \begin{proposition}[\cite{7}]
 \label{pro}
 Let $C\in\mathcal{M}_{m}$  such that $||  C||  <1$, then the matrix
$I-C$  is invertible and we have
\begin{equation}
||(I - C)^{-1}||  \leq \frac{1}{1 - ||C|| }.
 \end{equation}
\end{proposition}

To end this section, we give the following Theorem.

\begin{theorem}[\cite{11}] \label{gan}
 If the function $f(x)$  can be expanded in a~power series in
the circle $| x - x_{0}|<r$, namely
\begin{equation}
f(x) = \sum\limits_{\substack{p = 0}}^{+\infty}{\alpha_{p}(x -
x_{0})^{p}},
\end{equation}
then this expansion remains valid when the scalar argument $x$  is
replaced by a~matrix $A$  whose characteristic values lie within the
circle of convergence.
\end{theorem}

\section{Main results}
\subsection{The real case}

In mathematics, the error function, also called the Gauss error
 function, is a~special function of sigmoid shape that occurs in probability, statistics and partial differential equations describing diffusion. It is defined as
$$
\erf(x) = \frac{2}{\sqrt{\pi}}\displaystyle \int_{0}^{x} e^{-t^{2}} \, \mathrm{d}t, \hspace{0.5cm} |x|<\infty
$$
where the coefficient in front of the integral normalizes $\erf(+\infty) = 1$. A plot of $\erf(x)$  over
the range $-3\leq x \leq 3$  is shown as follows.
\begin{figure}
\begin{center}
\includegraphics[width = 5cm]{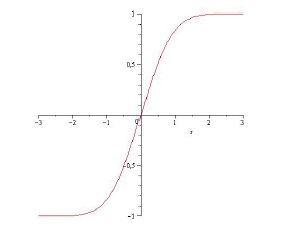}
\end{center}
\caption{plot of function error.}
\end{figure}

The power series expansion for the error function is given by
 \begin{equation}
 \erf(x) = \frac{2}{\sqrt{\pi}}\sum\limits_{\substack{n = 0}}^{+\infty}{\frac{(-1)^{n}}{(2n+1)\fact{n}}}x^{2n+1}.
 \label{sum1}
 \end{equation}
 Accordingly, we have
 \begin{equation}
 \erf(x) = \frac{2}{\sqrt{\pi}}\sum\limits_{\substack{n = 0}}^{+\infty}{\frac{(-1)^{n}x^{n+1}}{(2n+1)\fact{n}}}x^{n}.
 \end{equation}

 \begin{lemma}\label{lem2}
 Let $x$  be a~real number. Then the continued fraction expansion of the error function is
 \begin{equation}
\erf(x) = \left[0;\frac{(2/\sqrt{\pi})x}{1}, \frac{x^{2}}{3-x^{2}},
\frac{-(n-1)(2n-1)^{2}x^{2}}{(-1)^{n-1}(n(2n+1)-(2n-1)x^{2})}
\right]^{+\infty}_{n = 2}. \label{frac1}
\end{equation}
\end{lemma}

\begin{proof}
We use Lemma \ref{lemma} for the function
$$
g(x) = \sum\limits_{\substack{n = 0}}^{+\infty}{\frac{(-1)^{n}x^{n+1}}{(2n+1)\fact{n}}}x^{n},\quad c_{n} = \frac{(-1)^{n}x^{n+1}}{(2n+1)\fact{n}}.
$$
So, we have
$$
\frac{c_{0}}{1} = \frac{x}{1},\quad \frac{-c_{1}x}{c_{0}+c_{1}x} = \frac{\frac{x^{3}}{3}}{\frac{3x-x^{3}}{3}},\quad
\frac{-c_{0}c_{2}x}{c_{1}+c_{2}x} = \frac{\frac{-x^{5}}{10}}{\frac{-10x^{2}+3x^{4}}{30}}\,.
$$

For $n\geq 3$, we get
 \begin{align*}
c_{n-2}c_{n}x &= \frac{(-1)^{n-2}x^{n-2+1}}{(2(n-2)+1)\fact{(n-2)}}\cdot \frac{(-1)^{n}x^{n+1}}{(2n+1)\fact{n}} x\\
 &= \frac{x^{2n+1}}{n(n-1)(2n+1)(2n-3)(\fact{(n-2)})^{2}}\,.
\end{align*}
Furthermore, we have
\begin{align*}
c_{n-1}+c_{n}x &= \frac{(-1)^{n-1}x^{n-1+1}}{(2(n-1)+1)\fact{(n-1)}} + \frac{(-1)^{n}x^{n+1}}{(2n+1)\fact{n}}x\\
 &= \frac{(-1)^{n-1}x^{n}}{\fact{(n-1)}}(\frac{1}{2n-1}+\frac{-x^{2}}{n(2n+1)})\\
 &= \frac{(-1)^{n-1}(n(2n+1)-(2n-1)x^{2})}{(2n+1)(2n-1)\fact{n}}x^{n}\,.
\end{align*}
Then, we obtain
\begin{equation*}
\frac{-c_{n-2}.c_{n}x}{c_{n-1}+c_{n}x} =
\frac{\frac{-x^{2n+1}}{n(n-1)(2n+1)(2n-3)(\fact{(n-2)})^{2}}}{\frac{(-1)^{n-1}(n(2n+1)-(2n-1)x^{2})}{(2n+1)(2n-1)\fact{n}}x^{n}}.
\end{equation*}

Therefore, the continued fraction expansion of
$\erf(x) = \frac{2}{\sqrt{\pi}}g(x)$  is
\begin{align*}
\erf(x) &= \left[0;\frac{b_{n}}{a_{n}}\right]^{+\infty}_{n = 1} \\
 &= \left[0; \frac{(2/\sqrt{\pi})x}{1},
\frac{\frac{x^{3}}{3}}{\frac{3x-x^{3}}{3}},
\frac{\frac{-x^{5}}{10}}{\frac{-10x^{2}+3x^{4}}{30}},\frac{\frac{-x^{2n+1}}{n(n-1)(2n+1)(2n-3)(\fact{(n-2)})^{2}}}{\frac{(-1)^{n-1}(n(2n+1)-(2n-1)x^{2})}{(2n+1)(2n-1)\fact{n}}x^{n}}
\right]^{+\infty}_{n = 3}.
\end{align*}
Let us define the sequence
$(r_{n})_{n\geq1}$  by
$$
\begin{cases}
 r_{1} = 1\\
 r_{n} = \frac{(2n-1)(2n-3)(\fact{(n-1)})}{x^{n-1}}, \text{ for }  n\geq2.
 \end{cases}
$$
Then, we have 

$$
\begin{cases}
 \dfrac{r_{1}b_{1}}{r_{1}a_{1}} = \dfrac{(2/\sqrt{\pi})x}{1},\\
\dfrac{r_{1}r_{2}b_{2}}{r_{2}a_{2}} = \dfrac{x^{2}}{3-x^{2}}, \\
 \dfrac{r_{n}r_{n+1}b_{n+1}}{r_{n+1}a_{n+1}} = \dfrac{-(n-1)(2n-1)^{2}x^{2}}{(1-)^{n-1}(n(2n+1)-(2n-1)x^{2})} \text{ for } n\geq2.
 \end{cases}
$$

By applying the result of Lemma~\ref{lem} to the sequence
$(r_{n})_{n\geq1}$, we obtain
\begin{align*}
\erf(x) &= \left[0;\frac{(2/\sqrt{\pi})x}{1},
\frac{x^{2}}{3-x^{2}},\frac{-9x^{2}}{-10+3x^{2}},
\frac{-(n-1)(2n-1)^{2}x^{2}}{(-1)^{n-1}(n(2n+1)-(2n-1)x^{2})}
\right]^{+\infty}_{n = 3}\\
 &= \left[0;\frac{(2/\sqrt{\pi})x}{1},
\frac{x^{2}}{3-x^{2}},
\frac{-(n-1)(2n-1)^{2}x^{2}}{(-1)^{n-1}(n(2n+1)-(2n-1)x^{2})}
\right]^{+\infty}_{n = 2} \\
\end{align*}
 and the proof is complete.
\end{proof}

\subsection{The matrix case}
According to Theorem~\ref{gan}, we have
\begin{definition}
Let $A$  be a~matrix in $\mathcal{M}_{m}$. Then we define the error function
by the expression
\begin{equation}
\erf(A) = \frac{2}{\sqrt{\pi}}\sum\limits_{\substack{n = 0}}^{+\infty}{\frac{(-1)^{n}}{(2n+1)\fact{n}}}A^{2n+1}.
\label{sum2}
\end{equation}
\end{definition}
Now, we treat the matrix case,

\begin{theorem}\label{frac2}
Let $A$  be a~matrix in $\mathcal{M}_{m}$, such that $||  A||  = \alpha$, 
where $0<\alpha<\frac{1}{2}$. The continued fraction
$$
\left[0;\frac{(2/\sqrt{\pi})A}{I}, \frac{A^{2}}{3I-A^{2}},
\frac{-(n-1)(2n-1)^{2}A^{2}}{(-1)^{n-1}(n(2n+1)I-(2n-1)A^{2})}
\right]^{+\infty}_{n = 2}
$$
converges in $\mathcal{M}_{m}$. Furthermore, this continued fraction represents $\erf(A)$. So
$$
\erf(A) = \left[0;\frac{(2/\sqrt{\pi})A}{I}, \frac{A^{2}}{3I-A^{2}},
\frac{-(n-1)(2n-1)^{2}A^{2}}{(-1)^{n-1}(n(2n+1)I-(2n-1)A^{2})}
\right]^{+\infty}_{n = 2}.
$$
\begin{proof}
We study the convergence of the continued fraction $K(B_{k}/A_{k})$ 
with
$$
\begin{cases}
A_{1} = I, A_{2} = 3I-A^{2}, \\
B_{1} = (2/\sqrt{\pi})A, B_{2} = A^{2},
\end{cases}
$$
and for $k \geq 3$, we have:
$$
\begin{cases}
 A_{k} = (-1)^{k-2}((k-1)(2k-1)I-(2k-3)A^{2}),\\
B_{k} = -(k-2)(2k-3)^{2}A^{2},
\end{cases}
$$
we check that the conditions of Theorem~\ref{cov} are satisfied:
\begin{align*}
B_{2k-2} \ldots B_{2} &= \pm ((2k-2)-2)(2(2k-2)-3)^{2} \ldots (4-2)(2\cdot4-3)^{2} A^{2(k-1)} \\
&= \pm (2k-4)(4k-7)^{2} \ldots 50 A^{2(k-1)}, \\
A^{-1}_{2k-1} &= -((2k-2)(4k-3)I-(4k-5)A^{2})^{-1}
\end{align*}
and
\begin{align*}
B_{2k-1} {}&{} B_{2k-3} \ldots B_{1} \\
&= \pm(2k-1-2)(2(2k-1)-3)^{2} \ldots (3-2)(2\cdot3-3)^{2}\cdot
(2/\sqrt{\pi})A^{2k-1} \\
&= \pm(2/\sqrt{\pi})(2k-3)(4k-5)^{2}\cdot \ldots \cdot9 A^{2k-1}.
\end{align*}
Then, we have:
\newline
\resizebox{\textwidth}{!}{
\parbox{1.39\textwidth}{
\begin{align*}
\|{}&{}(B_{2k-2} \dots B_{2})^{-1}A^{-1}_{2k-1}B_{2k-1}B_{2k-3} \dots B_{1}\| \\
&{ }= \left\| \frac{1}{(2k-4)(4k-7)^{2} \ldots 50}A^{-2(k-1)} ((2k-2)(4k-3)I-(4k-5)A^{2})^{-1}(2/\sqrt{\pi})(2k-3)(4k-5)^{2} \dots 9 A^{2(k-1)+1}\right\| \\
&{ }\leq \frac{(2/\sqrt{\pi})(2k-3)(4k-5)^{2} \dots 9}{(2k-4)(4k-7)^{2} \dots 50} \| A^{-2(k-1)}((2k-2)(4k-3)I-(4k-5)A^{2})^{-1} A^{2(k-1)+1} \|.
\end{align*}
}
}
 
Now, the matrices
$((2k-2)(4k-3)I-(4k-5)A^{2})^{-1}$  and
$A^{-2(k-1)}$  commute, so the above inequality becomes 
\begin{multline*}
\left\| (B_{2k-2} \ldots B_{2})^{-1}A^{-1}_{2k-1}B_{2k-1}B_{2k-3} \ldots B_{1} \right\| \\
\leq\frac{(2/\sqrt{\pi})(2k-3)(4k-5)^{2} \ldots 9}{(2k-4)(4k-7)^{2} \ldots 50}
\left\| (I-\frac{(4k-5)}{(2k-2)(4k-3)}A^{2})^{-1} A \right\|.
\end{multline*}

By Proposition~\ref{pro} and the fact that 
$|| A||  <1/2$, we obtain
\[
\left\| (I-\frac{(4k-3)}{(2k-1)(4k-1)}A^{2})^{-1} \right\|
\leq \frac{1}{1 - \left\| \frac{(4k-5)}{(2k-2)(4k-3)}A^{2} \right\|} <1
\]
It implies that for all sufficiently large k, we get
$$
||B_{2k-2} \ldots B_{2})^{-1}A^{-1}_{2k-1}B_{2k-1}B_{2k-3} \ldots B_{1}||
\leq ||A||  = \alpha< 1/2.
$$
To prove the second inequality of Theorem~\ref{cov}, we have
\begin{multline*}
(B_{2k-1}B_{2k-3} \ldots B_{1})^{-1}A^{-1}_{2k}B_{2k} \ldots B_{2} \\
= \frac{(2k-2)(4k-3)^{2} \ldots
50}{(2/\sqrt{\pi})(2k-3)(4k-5)^{2} \ldots 9}
\frac{A^{2(k-1)+2}}{A^{2(k-1)+1}}
\left( I-\frac{4k-3}{(2k-1)(4k-1)}A^{2} \right)^{-1}
\end{multline*}
Again using the fact that the matrices
$\left(I-\frac{4k-3}{(2k-1)(4k-1)}A^{2}\right)^{-1}$  and $A^{2k-1}$  commute,
the Proposition~\ref{pro} and passing to the norm, we get
$$
|| (B_{2k-1}B_{2k-3} \ldots B_{1})^{-1}A^{-1}_{2k}B_{2k} \ldots B_{2}|| 
\leq ||  A||  = \alpha< 1/2
$$
which
completes the proof.
\end{proof}
\end{theorem}

\section{Numerical applications}
This section will provide some numerical data to illustrate the
preceding results. The focus will be on two cases:
\subsection{Real case:}

\begin{itemize}
\item The following table clarifies the differences between $\erf(x)$  and its first convergents when applying Lemma~\ref{lem2}.

\resizebox{.93\textwidth}{!}{
\begin{tabular}{|c|c|c|c|c|c|}
\hline
 x & $(\erf-F_{1})(x)$  & $(\erf-F_{2})(x)$  & $(\erf-F_{3})(x)$ & $(\erf-F_{4})(x)$  & $(\erf-F_{5})(x)/\erf(x)$  \\
 \hline
 0.005& -0.47015e-7 & 0.1e-11 & 0.1e-11 &0.1e-11& 0.1e-11 \\
 \hline
 0.05 & -0.00004698055 &0.3525e-7 & -0.2e-10 &0.1e-10& 0 \\
 \hline
 0.075 & -0.15841090e-3 & 0.26742e-6 & -0.35e-9 &-0.1e-10&-0.1e-10\\
 \hline
 0.1& -0.3750007e-3 & 0.11257e-5 & -0.27e-8 &0&0\\
 \hline
 0.15 & -0.12609036e-2 & 0.85229e-5 & -0.457e-7 &0.3e-9&0\\
 \hline
 0.2 & -0.29732442e-2 & 0.357669e-4 & -0.3412e-6 &0.27e-8& 0\\
 \hline
 0.25 & -0.57684016e-2 & 0.1085733e-3 & -0.16201e-5 &0.196e-7&-0.4e-9\\
 \hline
 0.3 & -0.98869906e-2 & 0.2684219e-3 & -0.57742e-5 &0.1014e-6&-0.12e-8\\
 \hline
 0.35 & -0.155506548e-1 & 0.5757641e-3 & -0.168817e-4 &0.4037e-6&-0.81e-8\\
\hline
 0.4 & -0.229593118e-1 & 0.11127771e-2 & -0.426832e-4 &0.13344e-5&-0.351e-7\\
 \hline
 0.45 & -0.322889054e-1 &0.19856118e-2 & -0.965652e-4 &0.38255e-5&-0.1275e-6\\
 \hline
\end{tabular}
}

We can clearly see that $F_{5}$  is approximately the
exact value of $\erf(x)$.
\item The following graphics illustrate the approximations of $\erf(x)$  in
terms of continued fractions.
\end{itemize}

\begin{figure}
\begin{center}
\includegraphics[width = 7cm]{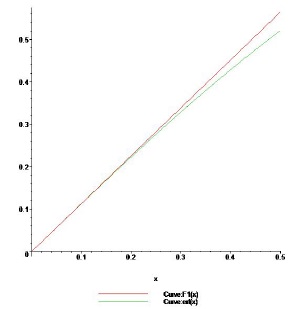}
\end{center}
\caption{Iteration 1.}
\end{figure}

\begin{figure}
\begin{center}
\includegraphics[width = 7cm]{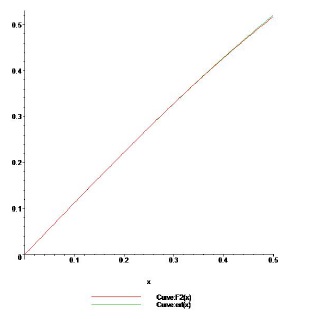}
\end{center}
\caption{Iteration 2.}
\end{figure}

\begin{figure}
\begin{center}
\includegraphics[width = 7cm]{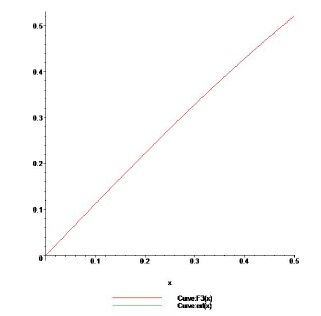}
\end{center}
\caption{Iteration 3.}
\end{figure}

After the first iteration, the convergence of
$(F_{n}(x))$  to $\erf(x)$  is very rapid. It is hard to distinguish
between the curve of the convergents and that of $\erf(x)$, for
$0\leq x\leq 0.5$.

\subsection{Matrix case:}

\begin{example}
Let $A$  be a~matrix such that
$$
A = 
\begin{pmatrix}
 \frac{1}{3} & \frac{1}{17} \\
 & \\
 \frac{-2}{23} & \frac{1}{11} \\
\end{pmatrix}
$$
The value of $\erf(A)$  is given by

$$
\erf(A) =
\begin{pmatrix}
 0.3640064111 & 0.06327099117 \\
 &\\
 -0.09353103045 & 0.1032532354 
\end{pmatrix}\,.
$$

 Using the expansion of Theorem~\ref{frac2}, we can
 obtain the following convergents of~$\erf(A):$ 
 $$
F_{1} = 
\begin{pmatrix}
 0.3636532973 & 0.06317676555 \\
 & \\
 -0.09339174038 & 0.1032884453 \\
\end{pmatrix}
\,.
$$

$$
F_{2} = 
\begin{pmatrix}
 0.3640145285 & 0.06327316896 \\
 & \\
 -0.09353424976 & 0.1032523777 \\
\end{pmatrix}
\,.
$$

$$
F_{3} = 
\begin{pmatrix}
 0.3640062588 & 0.06327095029 \\
 & \\
 -0.09353097002 & 0.1032532515 \\
\end{pmatrix}
\,.
$$
$$
F_{4} = 
\begin{pmatrix}
 0.3640064133 & 0.06327099181\\
 & \\
 -0.09353103134 & 0.1032532351 \\
\end{pmatrix}
\,.
$$

$$
F_{5} = 
\begin{pmatrix}
 0.3640064109 & 0.06327099114 \\
 & \\
 -0.09353103039 & 0.1032532353 \\
\end{pmatrix}
\,.
$$
\end{example}

\begin{example}

Let $A$  be a~matrix such that
$$
A = 
\begin{pmatrix}
 \frac{1}{15} & \frac{1}{9} & 0 \\
 & & \\
 0 & \frac{1}{20} & 0 \\
 & & \\
 \frac{1}{7} & 0 & \frac{1}{5} \\
\end{pmatrix}
\,.
$$
The value of $\erf(A)$  is given by
$$
\erf(A) = 
\begin{pmatrix}
 0.07511398139 & -0.1249466906 & 0 \\
 & & \\
 0 & 0.05637197780 & 0 \\
 & & \\
 0.1581306512 &0.0018630324 &0.2227025892 \\
\end{pmatrix}
\,.
$$

 We calculate $\erf(A)$  by using the expansion given in Theorem~\ref{frac2}.
 The first convergents are
 $$
F_{1} = 
\begin{pmatrix}
 0.07511383296&-0.1249459359 & 0\\
 & & \\
 0 & 0.05637194255 & 0 \\
 & & \\
 0.1580924886 & 0.001890582379 & 0.2226668223\\
\end{pmatrix}
\,.
$$
$$
F_{2} = 
\begin{pmatrix}
 0.07511398151 & -0.1249466915& 0 \\
 & & \\
 0 & 0.05637197782 &0 \\
 & & \\
 0.1581310167 & 0.001862762625 &0.2227029304 \\
\end{pmatrix}
\,.
$$

$$
F_{3} = 
\begin{pmatrix}
 0.07511398140 & -0.1249466906&0 \\
 & & \\
 0 & 0.05637197779 &0 \\
 & & \\
 0.1581306484 & 0.001863034561 & 0.2227025865\\
\end{pmatrix}
\,.
$$
$$
F_{4} = 
\begin{pmatrix}
 0.07511398140 & -0.1249466906 & 0 \\
 & & \\
 0 & 0.05637197779 & 0\\
 & & \\
 0.1581306513 & 0.001863032439 & 0.2227025892\\
\end{pmatrix}
\,.
$$

$$
F_{5} = 
\begin{pmatrix}
 0.07511398140 & -0.1249466906& 0 \\
 & & \\
 0 & 0.05637197779 & 0\\
 & & \\
 0.1581306512 & 0.001863032453 & 0.2227025892\\
\end{pmatrix}
\,.
$$

\end{example}

 \begin{example}
 Let $A$  be a~matrix such that
$$
A = 
\begin{pmatrix}
 0.1 & -0.02 & 0 & 0 & 0\\
 & & & & \\
 0 & 0.008 & 0 & 0 & 0 \\
 & & & & \\
 0.015 & -0.075 & 0.025 & -0.09 & 0 \\
 & & & & \\
 0.001 & 0 & 0 & 0.05 & 0 \\
 & & & & \\
 0.002 & 0 & 0 & 0.05 & 0.002\\
\end{pmatrix}
\,.
$$
The value of $\erf(A)$  is given by
\[
\resizebox{\textwidth}{!}{$
\erf(A) = 
\begin{pmatrix}
 0.1124632135 & -0.02248610288 & 0& 0.000007479569022 & 0.00001649965871\\
 & & & & \\
 -0.000002023944323 & 0.009027140204 & 0 & -0.00002255022425 & 0.02256695121 \\
 & & & & \\
 0.01685888359 &-0.08458907646 &0.02820360331 & -0.1013779889 & 0.00002199040709 \\
 & & & & \\
 0.001121818739 & 0.000001184002249 & 0& 0.05637197730& 1.495913804\, 10^{-7} \\
 & & & & \\
 0.002246257465 & 0.000002023944322 &0 & 0.05637002255& 0.002257054843 \\
\end{pmatrix}.
$}
\]

Now, let us apply the expansion of Theorem~\ref{frac2} to
obtain the following approximations of $\erf(A)$ 
\begin{gather*}
\resizebox{\textwidth}{!}{$
F_{1} = 
\begin{pmatrix}
 0.1128379167 &-0.02256758334 &0 & 0& 0 \\
 & & & & \\
 0 & 0.009027033336 &0 & 0& 0.02256758334 \\
 & & & & \\
 0.01692568750 & -0.08462843752& 0.02820947918& -0.1015541250 &0 \\
 & & & & \\
 0.001128379167 & 0 & 0 & 0.05641895835& 0 \\
 & & & & \\
 0.002256758334 & 0 & 0 &0.05641895835 & 0.002256758334 \\
\end{pmatrix},
$}
\\
\resizebox{\textwidth}{!}{$
F_{2} = 
\begin{pmatrix}
 0.1124620912 &-0.02248585859 & 0& 0.000007522527778&0.00001654956112 \\
 & & & & \\
 -0.000002031082501 & 0.009027141660 & 0& -0.00002256758334&0.02256695145 \\
 & & & & \\
 0.01685868999 & -0.08458902888 & 0.02820360220 &-0.1013778158 & 0.00002200339376\\
 & & & & \\
 0.001121796955 & 0.000001188559389& 0 & 0.05637194255& 1.504505556\, 10^{-7} \\
 & & & & \\
 0.002246223786 &0.000002031082501 & 0 & 0.05636998669& 0.002257056226 \\
\end{pmatrix},
$}
\\
\resizebox{\textwidth}{!}{$
F_{3} = 
\begin{pmatrix}
 0.1124632161 & -0.02248610346&0 &0.000007479459802 & 0.00001649953872\\
 && & & \\
 -0.000002023926862 &0.009027140199 & 0& -0.00002255019727 & 0.02256695120\\
 & & & & \\
 0.01685888405 & -0.08458907661 & 0.02820360330&-0.1013779889, &0.00002199031849] \\
 & & & & \\
 0.001121818792 & 0.000001183991006 & 0 &0.05637197734 & 1.495891960 \, 10^{-7} \\
 & & & & \\
 0.002246257546 & 0.000002023926862 & 0 & 0.05637002258& 0.002257054839 \\
\end{pmatrix},
$}
\\
\resizebox{\textwidth}{!}{$
F_{4} = 
\begin{pmatrix}
 0.1124632135 & -0.02248610288& 0& 0.000007479574705& 0.00001649965756 \\
 & & & & \\
 -0.000002023944490 & 0.009027140205 &0 &-0.00002255020712 &0.02256695120 \\
 & & & & \\
 0.01685888359 & -0.08458907650 & 0.02820360330& -0.1013779888&0.00002199033965 \\
 & & & & \\
0.001121818739 & 0.000001184002463& 0 & 0.05637197732 & 1.495914941\, 10^{-7} \\
 & & & & \\
0.002246257465 &0.000002023944490 & 0 & 0.05637002257 & 0.002257054844 \\
\end{pmatrix},
$}
\\
\resizebox{\textwidth}{!}{$
F_{5} = 
\begin{pmatrix}
 0.1124632135 &-0.02248610288 & 0 &0.000007479574479 &0.00001649965732] \\
 & & & & \\
 -0.000002023944456 & 0.009027140205 &0 & -0.00002255020711& 0.02256695120\\
& & & & \\
 0.01685888360 & -0.08458907650 &0.02820360330 & -0.1013779888 & 0.00002199033962\\
 & & & & \\
 0.001121818739 & 0.000001184002441& 0 &0.05637197732 & 1.495914896\, 10^{-7} \\
 & & & & \\
0.002246257465&0.000002023944456 &0 & 0.05637002257 & 0.002257054844 \\
\end{pmatrix}.
$}
\end{gather*}
\end{example}
In the examples above, we can clearly see that $F_{5}$  is
approximately the exact value of $\erf(A)$. This shows the
importance of the continued fractions approach.

\EditInfo{20 June, 2019}{28 January, 2020}{Karl Dilcher}

\end{paper}